\documentclass[journal]{IEEEtran}
\IEEEoverridecommandlockouts       

\usepackage{amsthm}
\usepackage{times} 
\usepackage{pdfsync} 
\usepackage{epsfig}
\usepackage{epstopdf}
\usepackage{amsmath} 
\usepackage{amssymb}
\usepackage{color} 
\usepackage{cite} 
\usepackage{array}
\usepackage{mathrsfs}
\usepackage[english]{babel} 
\usepackage{fix2col}

\newcommand{\field}[1]{\mathbb{#1}} 
\newcommand{\R}{\field{R}}
\newcommand{\N}{\field{N}}

\graphicspath{{Matlab/}}

\usepackage{algorithm}
\usepackage[noend]{algorithmic} 

\newtheorem{theorem}{Theorem}[section]

\newtheorem{proposition}[theorem]{Proposition}

\title{\LARGE \bf{Distributed Charging Control of Electric Vehicles\\ Using Online Learning}}

\author{Wann-Jiun Ma, Vijay Gupta, Ufuk Topcu
\thanks{W.-J. Ma and V. Gupta are with the Department of Electrical Engineering, University of Notre Dame, Notre Dame, IN 46556, USA. U. Topcu is with the Department of Electrical and Systems Engineering, University of Pennsylvania, Philadelphia, PA 19104, USA.  {\tt\small \{wma1, vgupta2\}@nd.edu,~\tt utopcu@seas.upenn.edu.} The work of W.-J. Ma and V. Gupta is supported in part by NSF grant $\#1239224$. The work of U. Topcu is supported in part by NSF grant $\#1312390$. 
 Initial results were presented at IEEE Conference on Decision and Control (CDC), 2014.      } 
}

\begin{document}
\maketitle

\begin{abstract}

We propose an algorithm for distributed charging control of electric vehicles (EVs) using online learning and online convex optimization. Many distributed charging control algorithms in the literature implicitly assume fast two-way communication between a distribution company and EV customers. This assumption is impractical at present and raises privacy and security concerns.
Our algorithm does not use this assumption; however, at the expense of slower convergence to the optimal solution. The proposed algorithm requires one-way communication, which is implemented through the distribution company publishing the pricing profiles of the previous days. We provide convergence results of the algorithm and illustrate the results through numerical examples.
\end{abstract}

\begin{IEEEkeywords}
Charging control, demand response, electric vehicles, online convex optimization, online learning, regret minimization
\end{IEEEkeywords}

\section{Introduction}\label{s1}

Demand response (DR) is an important functionality for the next generation power systems. It empowers the distribution company and its customers to decide collectively, but in a distributed manner, the best way to schedule energy usage. The reader is referred to \cite{cljl,gkgkzw,ri}, and references therein for the details of DR.   This paper focuses on the flexible load capability offered by electric vehicles (EVs) owned by residential customers. 


Large-scale integration of EVs imposes a significant burden on the grid. Particularly,  creation of new peaks, peak load amplification \cite{krw} and voltage deviations \cite{lsa} among other effects have been identified as major concerns. To cope with these issues, many algorithms have been proposed to schedule the charging of EVs, e.g., \cite{cf,chd,mch,mch2,shmv}. We are particularly interested in the algorithms such as those proposed in \cite{gtl,gwtcl,mgt,mch2}, which lead to the convergence of the total load profile to a desired one (for instance, a valley-filling profile) through appropriate price signals transmitted to the owners.

Many of the existing algorithms have analytical convergence guarantees and do not require the customers to share their charging constraints with the distribution company. On the other hand, they require a series of messages to be exchanged among the distribution company and the customers regarding possible price profiles and desired charging profiles in response.  As the available power supply and the customer requirements for charging their EVs change from day to day, these messages need to be exchanged daily to calculate the charging profiles. Since these message exchanges need to be completed before the EVs can begin charging, the algorithms, thus, implicitly assume the presence of a communication infrastructure and protocols that can support low-latency two-way communication between the distribution company and the EV customers. Such infrastructure and protocols have not yet been deployed extensively \cite{doe_comm}. Furthermore, the transmitted data carry information about the constraints faced by the individual customer and raises privacy and security concerns.  

Motivated by this issue, we propose an online learning and online convex optimization based distributed charging control algorithm. This algorithm requires only one-way communication from the distribution company to the customers. Furthermore, the communication carries information about the pricing profiles of \textit{previous} days.   


In our formulation, we model the distribution company and every EV customer as decision makers who wish to optimize their own utility functions. For the distribution company, the payoff is maximized if the total load profile over a day is valley-filling  \cite{gtl,gwtcl}. For the EV customer, the utility function  is maximized if the cost to charge the EV over a day is minimized. By designing a suitable pricing policy, the distribution company aims at ensuring that the charging profiles followed by the customers aggregate to a valley-filling profile.  
Our distributed charging control algorithm is based on an online learning and online convex optimization framework.  The only communication that occurs is when the distribution company notifies every customer of the pricing profiles incurred over the previous days.  The online learning framework has tremendous popularity in the online convex optimization and machine learning community (see e.g., \cite{bm,cl,hak,z,sss}, and references therein). We use a regret minimization algorithm \cite{rs} in the online learning framework. The regret minimization algorithm uses a \emph{regret} as the performance measure and provides an iterative way for every decision maker to update its policy such that, at convergence, the policy is optimal in a suitably defined sense. Informally, the regret minimization algorithm operates as follows. Consider a situation in which multiple decision makers need to design their own individual utility functions while satisfying coupled constraints. Furthermore, the utility functions of the decision makers may also be coupled through the actions of multiple decision makers and also the environment. The environment consists of the factors that the decision makers cannot control. The environment is uncertain and time-varying. Decisions are made repeatedly and the resulting payoffs are used to improve the decision policies used by the decision makers. Thus, in every iteration, every decision maker makes a decision given the current policy and the realization of the environment. Given the decisions of the decision makers and the realization of the environment, the various decision makers obtain certain payoffs. 

Our contributions are two-fold. First, we present a regret minimization based distributed algorithm for charging control of EVs that requires only one-way communication. Second, we allow heterogeneous  utility functions for the various EV owners, as would be the case when customers vary in the elasticity of shifting their loads in response to prices.  

Some relevant references that apply regret minimizetion to DR are \cite{kg,jtz,skg}. In \cite{kg}, real-time electricity pricing strategies for DR are designed using  regret minimization. However, the  focus of that work is on optimizing the utility function for the distribution company, and the customer behavior is  assumed to be such that the load change is linear in the price variation. {\color{black}  The objective of \cite{jtz} is to design pricing policies for the customers having price responsive loads. The exact demand function of the customer is assumed to be unknown to the pricing policy maker.} In \cite{jtz}, the distribution company is the only decision maker, whereas in our work, the distribution company and the EV customers are all decision makers.   In \cite{skg}, regret minimization is used to learn the  charging behavior of the EV customers. The price responsiveness for a community of customers is captured through a conditional random field model. The regret minimization algorithm is adopted to learn the parameters of the model.  

\par The remainder of the paper  is organized as follows.  
Section~\ref{s2}  presents the problem formulation. The basic framework and main results are presented in
Section~\ref{s3}. Section~\ref{s4} extends the basic framework to  some practical charging scenarios. 
Some numerical examples can be found in Section~\ref{s6}. Section~\ref{s7} concludes the paper.
     

\section{Problem Formulation}\label{s2}
We consider a scenario in which $N$ customers schedule the charging of their electric vehicles (EVs) daily. The charging needs to be completed over a day.    
Let there be $T$ time slots in a day and denote the set of these time slots by $\mathcal{T}:=\{1,...,T\}$.  Denote the set of EVs by $\mathcal{N}:=\{1,...,N\}$. Denote the base load on day $k$ by $D^{k}(t)\in \R, ~t\in \mathcal{T}$.   We assume that this base load is unknown to the EV customers and to the distribution company at the beginning of the $k$-th day when the charging schedules are fixed. Furthermore, the base load may vary from day to day. 
  Denote by $x^k_i(t)\in \R$ the charging rate of the $i$-th EV in the $t$-th time slot on the $k$-th day.   The charging profile of the $i$-th EV on the $k$-th day is denoted by a vector $x^k_i:=\big{(}x^k_i(1),x^k_i(2),...,x^k_i(T)\big{)}$. The aggregated charging profile of the EV customers is  described by a vector $x^k:=(x^k_1,...,x^k_N)$.   Let $x^{\mathrm{low}}_i(t)$ and $x^{\mathrm{up}}_i(t)$ denote the minimum and maximum charging rates for the $i$-th EV in the $t$-th time slot and $S_i$ denote the desired total charge for the $i$-th EV at the end of the $k$-th day; thus $S_{i}=\sum_{t\in\mathcal{T}}x_{i}^{k}(t)$, for every $k$. The total load  as seen by the distribution company is the sum of the base load and the charging rates adopted by the EVs.
 
The objective of the distribution company is to achieve a total load profile that is valley-filling while ensuring that both the base load and the EVs are supplied with the required amount of energy. The base load is inflexible, while the EV charging profile may be shaped as long as the desired amount of power is provided by the end of the day. The goal of the distribution company can be described as obtaining the aggregated charging profile $x^{k}$, for every $k\in\mathbb{N}_{>0}$,  that solves 
\begin{equation}\label{prob1}
\begin{array}{ll}
\underset{x^k}{\mathrm{minimize}}~~c_u^k(x^k)\\
\mathrm{subject~to}~~x_i^{\mathrm{low}}(t)\leq x^k_i(t)\leq x^{\mathrm{up}}_i(t),~ t\in \mathcal{T}, i\in \mathcal{N},\\
~~~~~~~~~~~~~~\sum_{t\in \mathcal{T}}x^k_i(t)=S_i,~ i\in \mathcal{N},
\end{array}
\end{equation}
 where the distribution company cost function $c_u^k$ is chosen as 
\begin{equation}\label{company_cost}
c_u^k(x^k):=\sum_{t\in \mathcal{T}}\big{(}D^k(t)+\sum_{i\in \mathcal{N}}x^k_i(t)\big{)} ^2. 
\end{equation}
The cost function (\ref{company_cost}) for the distribution company is also considered in \cite{gtl,gwtcl}.
By solving (\ref{prob1}), the distribution company can obtain a valley-filling total load profile while satisfying the charging requests of the customers. 


To incentivize the customers to choose charging profiles that in aggregate minimize the cost (\ref{company_cost}), the distribution company designs suitable pricing profiles for the power being supplied to the EVs. Every EV customer fixes the charging schedule at the beginning of the day based on the information about her own constraints and any information provided by the distribution company. 
A {\it{price-sensitive}} EV customer seeks to minimize the total cost of charging by suitably shaping her charging schedule. Thus, the optimization problem for each such customer $i,~i \in \mathcal{N}$ is to design a charging profile $x^k_i,~k\in\N_{>0}$ that solves 
\begin{equation}\label{prob2}
\begin{array}{ll}
\underset{x_i^k}{\mathrm{minimize}}~~c^{k}_i(x^k_i,x_{j\neq i}^k)\\
\mathrm{subject~to}~~x_i^{\mathrm{low}}(t)\leq x^k_i(t)\leq x^{\mathrm{up}}_i(t),~ t\in \mathcal{T},\\
~~~~~~~~~~~~~~\sum_{t\in \mathcal{T}}x^k_i(t)=S_i,
\end{array}
\end{equation}
where $c^k_i$ is a convex function in $x_i^k$ and is also a function of the other customers's charging profiles, where $x^k_{j\neq i}:=x^k_{j}$ for $j\in \mathcal{N},~j\neq i$. Since the pricing policy is possibly a function of the base load $D^k$ and other customer's charging profiles, $c_i^k$ inherits these features as well. 


The information flow is as follows. The distribution company monitors the total load and publishes the price profile for the previous day as realized according to a fixed and known pricing policy. The customers decide on the charging schedules for the next day with access to these pricing profiles for the previous days. No other communication occurs between the distribution company and the customers, or among the customers.



\section{Online Learning Framework}\label{s3}
We now adopt a  regret minimization framework to solve both  problem (\ref{prob1}) and problem (\ref{prob2}). The regret minimization framework is used for online learning and optimization and only requires one-way message exchanges between the distribution company and the customers. 


Let $L$ be a $\lambda$-strongly convex function with respect to a given norm $\|\cdot\|$. Let $D_L(\cdot,\cdot)$ denote the Bregman divergence \cite{rs} with respect to $L$. Let $\|\cdot\|_*$ denote the norm that is dual to $\|\cdot\|$. Let $\nabla L$ denote the gradient of $L$ and $\nabla L^{-1}$ denote the inverse mapping of $\nabla L$. For example, if $L$ is the squared Euclidean norm, i.e., $L(\cdot)=\|\cdot\|^2$, then the corresponding Bregman divergence $D_L(x,y)$ is equal to the squared Euclidean distance $\|x-y\|^2$, the inverse mapping $\nabla L^{-1}(x)$ is equal to $\frac{1}{2}x$, and the dual norm is the Euclidean norm. 
 
\textbf{Customer perspective:} For the $i$-th EV customer, the decision variable is her charging profile $x_i^k$ on the $k$-th day.
The set of feasible charging profiles is given as
\begin{equation}
\begin{split}
\mathcal{F}_i:=\Bigg{\{}x_i^k\in\R^T \mid x_i^{\mathrm{low}}(t)\leq x^k_i(t)\leq x^{\mathrm{up}}_i(t),\\t\in \mathcal{T},\sum_{t\in \mathcal{T}}x^k_i(t)=S_i\Bigg{\}}.
\end{split}
\end{equation}


We define $x_i^*$ as
\begin{equation}\label{sub_soln1}
x_i^*:=\mathrm{arg}\underset{x_i\in \mathcal{F}_i}{\text{min}}\sum_{k=1}^Kc^k_i(x_i).
\end{equation}

In the regret minimization framework. The notion of regret is used to measure the performance of an online algorithm \cite{z,sss}. For customer $i\in \mathcal{N}$, the customer regret after $K$ days $R_i$ is defined  as the difference between the cumulative cost function value of the charging profiles $x_i^k,~k=1,...,K$ generated by an online algorithm and the one generated by $x_i^*$, i.e., 
\begin{equation}\label{regret1} 
R_i(K,x_i^k):=\sum_{k=1}^Kc_i^k(x_i^k)-\underset{x_i\in \mathcal{F}_i}{\text{min}}\sum_{k=1}^Kc_i^k(x_i).
\end{equation}

As we mentioned in Section \ref{s2}, $c^k_i$ is possibly a function of the base load $D^k$ and other customer's charging profiles $x_{j\neq i},~j\in \mathcal{N},~j\neq i$. Since the base load and the charging profiles may change from day to day, $c_i^k$ may not remain the same from day to day a well. While $x_i^*$ remains the same from one day to the next, $x_i^*$ is a suboptimal solution of (\ref{prob2}). The regret $R_i$ measures the difference between the performance of the charging profile generated by an online algorithm and the performance that is obtained by the suboptimal charging profile $x_i^*$. Notice that the suboptimal charging profile $x_i^*$ can only be calculated in hindsight after $K$ days have elapsed.

We adopt the optimistic mirror descent (OMD) algorithm  \cite{rs} to generate the charging profile update which minimizes the regret (\ref{regret1}). On each day, the regret minimization algorithm generates the charging profile update without knowing the current objective function (and its gradient). Specifically, the OMD algorithm iteratively applies the updates
\begin{equation}\label{omd} 
\begin{split}
h_i^{k+1}&=\nabla L_i^{-1}\big{(}\nabla L_i(h_i^k)-\eta_i \nabla c_i^k(x_i^k)\big{)}\\
x_i^{k+1}&=\underset{x_i\in \mathcal{F}_i}{\text{argmin}}~ \eta_i x_i^TM_i^{k+1}+D_{L_i}(x_i,h_i^{k+1}),
\end{split} 
\end{equation}
where $\eta_i \in \R $ is an algorithm parameter, $h_i^k$ is an intermediate update of the charging profile. For easy of presentation, for the vector $h_i^k\in\R^{T}$, $L_i(h_i^k)$ is set to $L_i(h_i^k)=\frac{\|h_i^k\|^2}{2}$. The $i$-th customer may have a prediction $M_i^k$ for the gradient of the cost function $c_i^k$. For example, on day $k$, for customer $i\in \mathcal{N}$, one possible option of the prediction $M^k_i$ is the average of the gradients of the cost functions for the previous days, namely, $M_i^k(x_i^k)=\frac{1}{k-1}\sum_{\tilde{k}=1,...k-1}\nabla c_i^{\tilde{k}}(x_i^{\tilde{k}})$. 
Intuitively, the iteration (\ref{omd}) updates the charging profile toward the negative gradient direction and projects it onto the set of feasible charging profiles. 

 
As we mentioned in Section \ref{s2}, the customer and the distribution company have different objectives and hence different regrets. We now switch to the distribution company perspective and define the regret minimization framework for the company.

\textbf{Company perspective:} For the distribution company, the decision variable is the aggregated charging profile $x^k$ on the $k$-th day. The set of the aggregated feasible charging profiles is  denoted  by  $\mathcal{F}:=\mathcal{F}_1\times \mathcal{F}_2\times...\times\mathcal{F}_N$. 
%
The distribution company's regret after $K$ days is given by
\begin{equation}\label{regret2} 
R_u(K,x^k):=\sum_{k=1}^Kc_u^k(x^k)-\underset{x\in \mathcal{F}}{\text{min}}\sum_{k=1}^Kc_u^k(x).
\end{equation}
We define $x^*$ as
\begin{equation}\label{sub_soln2}
x^*:=\mathrm{arg}\underset{x\in \mathcal{F}}{\text{min}}\sum_{k=1}^Kc^k_u(x).
\end{equation}
  
The OMD algorithm generates the charging profile update which minimizes the regret  (\ref{regret2}) as
\begin{equation}\label{omdu} 
\begin{split}
h_u^{k+1}&=\nabla L_u^{-1}\big{(}\nabla L_u(h_u^k)-\eta_u \nabla c_u^k(x^k)\big{)},\\
x^{k+1}&=\underset{x\in \mathcal{F}}{\text{argmin}}~ \eta_u x^TM_u^{k+1}+D_{L_u}(x,h_u^{k+1}),
\end{split} 
\end{equation}
where $\eta_u \in \R $ is an algorithm parameter, $h_u^k$ is an intermediate update of the aggregated  charging profile, and $M^k_u$ is the prediction of the gradient of the cost function $\nabla c_u^k$. For easy of presentation, for the vector $h_u^k\in\R^{NT}$, $L_i(h_u^k)$ is set to $L_u(h_u^k)=\frac{\|h_u^k\|^2}{2}$. As an example, the prediction $M^k_u$ can be chosen to be the average of the gradients of the cost functions for the previous days, namely, $M_u^k(x^k)=\frac{1}{k-1}\sum_{\tilde{k}=1,...k-1}\nabla c_u^{\tilde{k}}(x^{\tilde{k}})$. 

\subsection{Convergence Results}
The following results summarize the convergence of the charging profile updates generated by the OMD algorithm. All proofs can be found in the Appendix.
\begin{proposition}\label{prop1}
(Convergence of regret): For every $x_i^{*}\in \mathcal{F}_i$, the iteration (\ref{omd}) converges in the sense that
\begin{equation}\label{upbd1}
 R_i(K,x_i^k)\leq  \frac{1}{\eta_i}P_i+\frac{\eta_i}{2}\sum_{k=1}^{K}\|\nabla c_i^k(x_i^k)-M_i^k\|_{*}^2,
 \end{equation}
 where
 \begin{equation}
 \begin{split}
 &R_i(K,x_i^k):=\sum_{k=1}^Kc_i^k(x_i^k)-\sum_{k=1}^Kc_i^k(x_i^*),\\
   &P_i:= \underset{x_i\in \mathcal{F}_i}{\emph{\text{max}}}~L_i(x_i)-\underset{x_i\in \mathcal{F}_i}{\emph{\text{min}}}~L_i(x_i).
   \end{split}
  \end{equation}
  In  particular, if $\eta_i$ is chosen as $O(1/{\sqrt{K}})$, then the average regret, i.e., $R_i(K)/K$, converges to zero as $K\rightarrow \infty$.\\
Similarly, for every $x^{*}\in \mathcal{F}$, the iteration (\ref{omdu}) converges in the sense that
\begin{equation}\label{upbd2}
 R_u(K,x^k)\leq  \frac{1}{\eta_u}P_u+\frac{\eta_u}{2}\sum_{k=1}^{K}\|\nabla c_u^k(x^k)-M_u^k\|_{*}^2,
 \end{equation}
 where
 \begin{equation}
 \begin{split}
 R_u(K,x^k):=\sum_{k=1}^Kc_u^k(x^k)-\sum_{k=1}^Kc_u^k(x^{*}),\\
   P_u:= \underset{x\in \mathcal{F}}{\emph{\text{max}}}~L_u(x)-\underset{x\in \mathcal{F}}{\emph{\text{min}}}~L_u(x).
   \end{split}
  \end{equation}
  In  particular, if $\eta_u$ is chosen as $O(1/{\sqrt{K}})$, then the average regret, i.e., $R_u(K)/K$, converges to zero as $K\rightarrow \infty$.
\end{proposition}


Any online algorithm yields a sublinear regret bound as in (\ref{upbd1}) and (\ref{upbd2}) is called a no regret algorithm \cite{cl,sss}. The results (\ref{upbd1}) and (\ref{upbd2}) guarantee that, as the number of days increases, the average performance of the charging profiles generated by the OMD algorithm approaches the performance that is obtained by the charging profiles $x^*$ and $x_i^*,~i\in \mathcal{N}$, respectively. Note that while the charging profiles $x^*$ and $x_i^*,~i\in \mathcal{N}$ may not solve the optimization problems (\ref{prob1}) and (\ref{prob2}), respectively, these solutions are optimal for the related problems (\ref{sub_soln1}) and (\ref{sub_soln2}), respectively. In Section \ref{day_varying}, we compare the performance of the charging profile generated by the OMD algorithm with the performance that is obtained by the solutions of problems (\ref{prob1}) and (\ref{prob2}). However, in that case, the convergence guarantees that are obtained are weaker.


\subsection{Design of the Pricing Function}

There are no guarantees that the solutions $x_i^*,~i\in\mathcal{N}$ of the problem (\ref{sub_soln1}) can solve the problem (\ref{sub_soln2}). In fact, unless the pricing function $c_i^k$ is carefully designed, these solutions will not be the same since the objectives of the distribution company and the EV customers are different. After some algebraic manipulation of the updates (\ref{omd}) and (\ref{omdu}), we observe that the natural choice of $c_i^k$  as 
\begin{equation}\label{cost_price}
 c_i^k(x_i^k)=\bigg(\sum_{j=1}^{N}x_j^k+D^k\bigg)^{T}x^k_i
 \end{equation}
does not lead to the charging profiles $(x_1^*,...,x_N^*)$ that reduce the regret of the distribution company to zero.  

 
We now propose a choice of $c_i^k$ to ensure that when each customer minimizes her regret, the aggregated charging profile minimizes the distribution company's regret.
\begin{proposition}\label{align}
If $c_i^k$ is chosen as
\begin{equation}\label{cust2}
c_i^k(x_i^k)=\bigg(\frac{1}{2}x_i^k+\sum_{j\neq i}^{N}x_j^k+D^k\bigg)^{T}x^k_i,~\quad i\in \mathcal{N},
\end{equation} 
the customers adopt the iteration (\ref{omd}), and $\eta_u=\frac{1}{2}\eta_i$, then the average regret of the distribution company as defined in (\ref{regret2})   converges to zero as the total number of days goes to infinity.
\end{proposition}

To update the charging profile on day $k$, the $i$-th customer needs to know $2x_i^{k-1}+\sum_{j\neq i}x_j^{k-1}+D^{k-1}$ or $\sum_{j}x_j^{k-1}+D^{k-1}$ depending on whether the pricing function (\ref{cost_price}) or (\ref{cust2}) is adopted. 
 The distribution company can simply publish the total load information for the previous day. The customers do not need to have full knowledge about how their consumption will map to a corresponding expenditure.    

\section{Extensions}\label{s4}  
The basic framework presented above can be extended in various directions. They include considering a different definition of regret and incorporating customers who vary in the elasticity of shifting their EV charging load in response to price. 

To ensure that the distribution company's average regret has a convergent behavior and the regret is of the order $O(\sqrt{K})$, in the following discussion we assume that the distribution company selects the cost function (\ref{cust2}) for each EV customer and sets $\eta_u=\frac{1}{2}\eta_i,~i\in \mathcal{N}$.

\subsection{Regret with Respect to the Optimal Charging Profiles}\label{day_varying}
The regrets defined in (\ref{regret1}) and (\ref{regret2}) measure the difference between the performance of the charging profiles generated by our algorithm and the performance that is obtained by the charging profiles $x^*$ and $x_i^*,~i\in \mathcal{N}$ that are the solutions of the related optimization problems (\ref{sub_soln1}) and (\ref{sub_soln2}), respectively. We can instead consider the original optimization problems (\ref{prob1}) and (\ref{prob2}) to define \emph{tracking regret} after $K$ days as
\begin{equation}\label{tracking_regret} 
R_u^{\text{tracking}}(K,x^k):=\sum_{k=1}^Kc_u^k(x^k)-\underset{x^k\in \mathcal{F},\forall k\in\mathcal{K}}{\text{min}}\sum_{k=1}^Kc_u^k(x^k),
\end{equation}
where $\mathcal{K}:=\{1,...,K\}$. We define the set $\{x^{k*},~ k\in \mathcal{K}\}$ as 
\begin{equation}
\Bigg{\{}x^{k*}\in\R^{NT},~k\in \mathcal{K}\mid x^{k*}=\mathrm{arg}\underset{x^k\in \mathcal{F},\\\\ \forall k \in \mathcal{K}}{\text{min}}\sum_{k=1}^Kc_u^k(x^k)\Bigg{\}}. 
\end{equation}
   
This notion of tracking regret characterizes the difference between the cumulative cost of the charging profiles generated by our algorithm and the cumulative cost of executing the optimal charging profiles that can be calculated only in hindsight.   For comparison, we refer to the the regrets (\ref{regret1}) and (\ref{regret2}) as \emph{static regrets}.  

%

\begin{theorem}\label{thm:tracking}
 For every $x^{k*}\in \mathcal{F}$, the OMD algorithm yields that,
\begin{equation}\label{tracking_upbd}
\begin{split}
&R_u^{\emph{\text{tracking}}}(K,x^k)\\
&\leq  \frac{1}{\eta_u}\bigg{[}L_u(h_u^{K+1})-L_u(h_u^1)\bigg{]}\\
&+\frac{1}{\eta_u}\bigg{[}\nabla L_u(h_u^{K+1})^T(x^{K+1*}-h_u^{K+1})\\
&-\nabla L_u(h_u^{1})^T(x^{1*}-h_u^{1})\bigg{]}\\
&+\frac{1}{\eta_u}\underset{k\in\mathcal{K}}{\emph{\text{max}}}\|\nabla L_u(h_u^k)\| \sum_{k=1}^{K}\|x^{k*}-x^{k+1*}\|\\
&+\frac{\eta_u}{2}\sum_{k=1}^{K}\|\nabla c_u^k(x^k)-M_u^k\|_{*}^2.
\end{split}
\end{equation}
In particular, if $\eta_u=O(1/\sqrt{K})$, then the tracking regret is of the order $O\big{(}\sqrt{K}[1+\sum_{k=1}^K\|x^{k*}-x^{k+1*}\|]\big{)}$. 
\end{theorem}
A comparison of the regret bounds (\ref{upbd2}) and (\ref{tracking_upbd}) is of interest. If $\eta_u=O(\sqrt{K})$, in  (\ref{upbd2}), the first term $\frac{1}{\eta_u}P_u$ and the second term $\frac{\eta_u}{2}\sum_{k=1}^{K}\|\nabla c_u^k(x^k)-M_u^k\|_{*}^2$ are of order $O(\sqrt{K})$. In (\ref{tracking_upbd}), the first three terms  measure the difference between the initial iterate $h_u^1$ and the final iterate $h_u^{K+1}$, the difference between the final iterate $h_u^{K+1}$ and the optimal solution $x^{K+1*}$, and the difference between the iterate $h_u^{1}$ and the optimal solution $x^{1*}$. The first three terms are of order $O(\sqrt{K})$. The last term in (\ref{tracking_upbd}), $\frac{\eta_u}{2}\sum_{k=1}^{K}\|\nabla c_u^k(x^k)-M_u^k\|_{*}^2$, also appears in (\ref{upbd2}) and is of order $O(\sqrt{K})$. The fourth term increases as $K$ increases and is of order  $O\big{(}\sqrt{K}\sum_{k=1}^K\|x^{k*}-x^{k+1*}\|\big{)}$. Because of the presence of the fourth term, the regret bound in (\ref{tracking_upbd}) increases as the variation of the optimal sequence of decisions $\sum_{k=1}^K\|x^{k*}-x^{k+1*}\|$ increases.
 If the optimal solution remains the same from one day to the next, then $\sum_{k=1}^K\|x^{k*}-x^{k+1*}\|=0$ and the  tracking regret $R_u^{\mathrm{tracking}}$ is of the order $O(\sqrt{K})$. On the other hand, if the optimal solution varies significantly from one day to the next, then the tracking regret $R_u^{\mathrm{tracking}}$ will be of the order $O\big{(}\sqrt{K}[1+\sum_{k=1}^K\|x^{k*}-x^{k+1*}\|]\big{)}$, and the average tracking regret will not necessarily converge to zero.   
 
If the distribution company has perfect prediction of the gradient of the cost function, i.e., $\nabla c_u^k(x^k)=M_u^k$ for all $k$, then the last term in (\ref{tracking_upbd}) vanishes and the distribution company can set $\eta\rightarrow \infty$ to ensure that the regret bound is zero, i.e.,  $R_u^{\mathrm{tracking}}\leq0$. It indicates that the cumulative cost function value of the charging profiles $x^k,~k=1,...,K$ generated by the online algorithm (\ref{omdu}) and the one generated by the elements in the set $\{x^{k*}\in\R^{NT},~k\in \mathcal{K}\}$ are identical. Note that the elements in the set $\{x^{k*}\in\R^{NT},~k\in \mathcal{K}\}$ solve problem (\ref{prob1}).

\subsection{Presence of Inelastic Customers}\label{4a}
The discussion so far assumed that all customers were rational in the sense that they wanted to choose their charging profile to solve (\ref{prob2}). Furthermore, they were elastic in scheduling their charging (within the constraints pre-specified by $x_i^{\mathrm{low}}(t)$, $x_i^{\mathrm{up}}(t)$, and $S_i$). We now assume that some customers are either irrational or inelastic and they do not optimize their schedules to solve (\ref{prob2}).
Suppose that $N_l$ out of $N$ customers are inelastic. Denote the set of inelastic customers by $\mathcal{N}_l$.  For every inelastic customer $i \in \mathcal{N}_l$, we assume that her charging profile remains the same from day to day and is not updated to minimize (\ref{prob2}). Equivalently, for the inelastic customers, the cost function $c^k_i$ can be  selected as $c_i^k\equiv r$, for all $k$, where $r$ is an arbitrary constant.  


Since the inelastic customers do not carry out any predictions, we set all customers' predictions to zeros, i.e., $M_i^k=0,~i\in \mathcal{N},~k\in \N_{>0}$. The update of the charging profile for inelastic  customer $i\in \mathcal{N}_l$ can thus be written as
\begin{equation} 
\begin{split}
h_i^{k+1}&=\nabla L_i^{-1}\Big{(}\nabla L_i(h_i^k)-\eta_i\big{(}\nabla c_i^k(x^k)+\epsilon_i^k\big{)}\Big{)},\\
x_i^{k+1}&=\underset{x_i\in \mathcal{F}_i}{\text{argmin}}~D_{L_i}(x_i,h_i^{k+1}),
\end{split} 
\end{equation}
where $\epsilon_i^k$ is an error term. For instance, for the cost function (\ref{cust2}), the error $\epsilon_i^k$ is equal to the total load  $\sum_{j}x_j^{k}+D^{k}$ on the $k$-th day. This error term quantifies the inconsistency between the updates as desired by the distribution company for each customer to execute and the inelastic customer's behavior. Denote by $\epsilon^k$ the aggregate error, i.e., $\epsilon^k:=(\epsilon_1^k,...,\epsilon_N^k)$.    


Due to the presence of the inelastic customers, the ability of the aggregated solution to be valley-filling and hence to minimize the cost function in problem (\ref{prob1}) is decreased. The performance loss is given in the following result.
\begin{theorem}\label{lazy}
Consider that there are $N$ EV customers out of which $N_l$ are inelastic customers. If $\eta_u=\frac{1}{2}\eta_i,~i\in\mathcal{N}$, then for every $x^{*}\in \mathcal{F}$, the regret of the distribution company  is bounded as
\begin{equation}\label{lazy_bd}
\begin{split}
 R_u(K,x^k)&\leq  \frac{1}{\eta_u}P_u+\frac{\eta_u}{2}\sum_{k=1}^{K}\big{\|}\big{(}\nabla c_u^k(x^k)+\epsilon^k\big{)}\big{\|}_{*}^2\\
 &+K\sum_{i\in\mathcal{N}_l}\|F_i\|\|\epsilon_i\|,
 \end{split}
 \end{equation}
      where
 \begin{equation}\label{pu}
  \begin{split}
    &P_u:= \underset{x\in \mathcal{F}}{\emph{\text{max}}}~L_u(x)-\underset{x\in \mathcal{F}}{\emph{\text{min}}}~L_u(x),\\
  &\|F_i\|:=\underset{x,y\in\mathcal{F}_i}{\text{\emph{max}}}\|x-y\|,\quad\|\epsilon_i\|:=\underset{k}{\text{\emph{max}}}\| \epsilon_i^k \|,~i\in\mathcal{N}_l.
     \end{split}
  \end{equation}
\end{theorem}

In (\ref{lazy_bd}), the average regret converges to a constant, i.e., $\lim_{K \to \infty} R_u(K)/K=\sum_{i\in\mathcal{N}_l}\|F_i\|\|\epsilon_i\|$.  The size of this constant depends on the error terms $\epsilon_i,~i\in\mathcal{N}_l$ and the charging constraints of the inelastic customers. The result explicitly quantifies the deviation from the desired performance in terms of the asymptotic average regret of the distribution company. To obtain further insight into the effect of the inelastic customers, we proceed as follows.

Note that the bound (\ref{lazy_bd}) depends on the size of the term $\|\epsilon_i\|$. For inelastic customer  $i\in \mathcal{N}_l$, the error term $\|\epsilon_i\|$ can be bounded as 
\begin{equation}
\begin{split}
&\|\epsilon_i\| = \underset{k}{\text{max}}\Bigg{\|}\Bigg(\sum_{i\in\mathcal{N}}x_i^k+D^k\Bigg)\Bigg{\|}\leq \underset{k}{\text{max}}\Bigg{(}\sum_{i\in\mathcal{N}}\|x_i^k\|+\|D^k\|\Bigg{)}\\
& \leq\sum_{i\in\mathcal{N}}\|x_i^{\mathrm{up}}\|+\underset{k}{\text{max}}\|D^k\|.
\end{split}
\end{equation}
Assume that the Euclidean ball $B$ is the smallest Euclidean ball containing the set of feasible charging schedules $\mathcal{F}_i$ and $r$ is the radius of the ball $B$. Then, the term $\|F_i\|$ can be further bounded by $2r$.  
Since the set $\mathcal{F}_i$ is a polytope, the algorithm in \cite{kmy} can be adopted to compute this bound efficiently. 

Notice that if $N_l=0$, (\ref{lazy_bd}) boils down to (\ref{upbd2}) (without the prediction, i.e., $M_u^k=0$ for all $k=1,...,K$). However, it becomes difficult to ensure the feasibility of the problem (\ref{prob1}). To improve the ability of the aggregated solution to be valley-filling  in the presence of inelastic customers, the distribution company can consider using some loads that can be controlled completely. We now discuss this option.   

\subsection{Controllable Customers}\label{control}

On the other end of the spectrum from customers that are completely inelastic are customers that are under the complete control of the distribution company. The customers under the complete control adopt the charging profiles assigned by the company. The charging constraints of these customers are also known and controlled  by the distribution company.   In practice, the distribution company can offer a contract to a subset of customers offering a special price to be such a controllable load. This contract-based direct load control has been implemented by many distribution
companies, e.g., \cite{con,wis}. 

To introduce such directly controlled customers in our formulation, we allow the charging constraints of the  controllable customers to be relaxed (thus enlarging the set of feasible charging profiles for the controllable customers). 
Denote the set of controllable customers by $\mathcal{N}_c$.  
On day $k$, for the controllable customer $i\in \mathcal{N}_c$, the set of the relaxed feasible charging profiles $\mathcal{\tilde{F}}_i$ is defined as
\begin{equation}
\begin{split}
\mathcal{\tilde{F}}_i:=\Bigg{\{}\tilde{x}_i^k\in\R^T\mid \tilde{x}_i^{\mathrm{low}}(t)\leq \tilde{x}^k_i(t)\leq \tilde{x}^{\mathrm{up}}_i(t),\\ 
t\in \mathcal{T},~\sum_{t\in \mathcal{T}}a_i\tilde{x}^k_i(t)=\tilde{S}_i\Bigg{\}},
\end{split}
\end{equation}
where $a_i=\{0,1\}$, $\tilde{x}_i^k$ is the feasible charging profile in the set $\tilde{\mathcal{F}_i}$, and $\tilde{x}_i^{\mathrm{low}}$, $\tilde{x}_i^{\mathrm{up}}$, $\tilde{S}_i$, the relaxed minimum charging rate, the relaxed maximum charging rate, and the relaxed total charge, respectively. The aggregated relaxed charging profile of the EV customers is  described by a vector $\tilde{x}^k$ which consists of all elements in the set
\begin{equation}
\Big\{\tilde{x}^k_i\in\R^T\mid i\in \mathcal{N}_c\Big\}\bigcup \Big\{x^k_i\in\R^T\mid i\in \mathcal{N}_l,~i\in \mathcal{N}_p\Big\},
\end{equation}
where $\mathcal{N}_p$ denotes the set of the price-sensitive customers.

For controllable customer $i\in\mathcal{N}_c$, we define $\tilde{x}_i^*$ as   
\begin{equation}
\tilde{x}_i^*:=\mathrm{arg}\underset{{\tilde{x}_i\in \mathcal{\tilde{F}}_i}}{\text{min}}\sum_{k=1}^Kc^k_i(\tilde{x}_i).
\end{equation}
We use $\tilde{\mathcal{F}}$ to denote the set of the feasible aggregated charging profiles including the controllable customers. We define  $\tilde{x}^*$ as  
\begin{equation}
\tilde{x}^*:=\mathrm{arg}\underset{\tilde{x}\in \tilde{\mathcal{F}}}{\text{min}}\sum_{k=1}^Kc^k_u(\tilde{x}).
\end{equation}

There are different ways to relax the set of feasible charging profiles for the controllable customers. For example, if we select $a_i=0$ and $\tilde{S}_i=0$, then the equality constraint is removed from the set of feasible charging profiles for the $i$-th customer, $i \in \mathcal{N}_c$, namely, the $i$-th controllable customer removes her total charging sum requirement daily. We can also relax the charging deadline by adjusting $\tilde{x}^{\mathrm{low}}_i(t)$ and $\tilde{x}^{\mathrm{up}}_i(t)$.

By enlarging the set of the feasible charging profiles, for controllable customer $i\in \mathcal{N}_c$, the cumulative cost function value of the iterates that are generated by the iteration (\ref{omd}) is a lower bound for the cumulative cost function value  of the iterates that are obtained by using  the same iteration (\ref{omd}) without relaxation.


We now propose Algorithm 1 with the above mentioned set of relaxed feasible charging schedules $\tilde{\mathcal{F}}_i,~i\in \mathcal{N}_c$ to compensate for the performance loss due to the presence of inelastic customers. We consider that the controllable customers relax their charging constraints for the final $J$ days. 

\begin{algorithm}[H]
\caption{}
\begin{algorithmic}[1]
\REQUIRE The distribution company knows $\mathcal{N}_p$, $\mathcal{N}_l$, $\mathcal{N}_c$, $\tilde{\mathcal{F}}_i,~i \in \mathcal{N}_c$, $J$, and sets $\eta_u=1/(2\sqrt{K})$. Each price-sensitive customer $i\in\mathcal{N}_p$ knows $\mathcal{F}_i$, and sets $\eta_i=1/\sqrt{K}$. Each inelastic customer $i\in\mathcal{N}_l$ knows $\mathcal{F}_i$. Each controllable customer $i\in\mathcal{N}_c$ knows $\mathcal{F}_i$, $\tilde{\mathcal{F}}_i$, $J$,  and sets $\eta_i=1/\sqrt{K}$
\STATE Initialization: $k \leftarrow 1$, $x^k_i(t) \leftarrow S_i/T$, for all $i\in \mathcal{N},~ t\in\mathcal{T}$
 
\FOR{$k=1$ to $K$}
\STATE At the end of the $k$-th day, the distribution company gathers the charging profiles $x_i^k,~i\in \mathcal{N}$, computes the prices $p_i^k(x^k)=(D^k+\sum_{i\in \mathcal{N}}x_i^k),~i\in \mathcal{N}$,  notifies the prices to the price-sensitive customers and the controllable customers
\STATE At the end of the $k$-th day, every price-sensitive customer $i \in \mathcal{N}_p$ updates the charging profile  by
 \begin{equation*} 
\begin{split}
h_i^{k+1}&=h_i^k-\eta_ip_i^k,\\
x_i^{k+1}&=\underset{x_i\in \mathcal{F}_i}{\text{argmin}}~\|x_i-h_i^{k+1}\|^2.
\end{split} 
  \end{equation*}
  \STATE At the end of the $k$-th day,   
every inelastic customer $i \in \mathcal{N}_l$ updates the charging profile  by
 \begin{equation*} 
x_i^{k+1}= x_i^k.
  \end{equation*}
  \STATE At the end of the $k$-th day, 
\IF{$1 \leq k \leq (K-J)$}
\STATE every controllable customer  $i \in \mathcal{N}_c$  updates the charging profile  by
  \begin{equation*} 
\begin{split}
h_i^{k+1}&=h_i^k-\eta_ip_i^k,\\
x_i^{k+1}&=\underset{x_i\in \mathcal{F}_i}{\text{argmin}}~\|x_i-h_i^{k+1}\|^2. 
\end{split} 
  \end{equation*}
\ELSE
\STATE every controllable customer  $i \in \mathcal{N}_c$ updates the charging profile  by
   \begin{equation*} 
\begin{split}
h_i^{k+1}&=h_i^k-\eta_ip_i^k,\\
\tilde{x}_i^{k+1}&=\underset{\tilde{x}_i\in \tilde{\mathcal{F}}_i}{\text{argmin}}~\|\tilde{x}_i-h_i^{k+1}\|^2. 
\end{split} 
  \end{equation*}
\ENDIF
 \ENDFOR
\end{algorithmic}
\end{algorithm} 

We have the following results to verify that Algorithm 1 is a no regret algorithm even with the presence of the inelastic customers.   

\begin{theorem}\label{lazy2}
Consider that there are $N$ EV customers out of which $N_l$ are lazy and $N_c$ are controllable. Assume that the controllable customers relax their charging constraints for the final $J$ days, $\eta_u=\frac{1}{2}\eta_i$, for all $i\in \mathcal{N}$, and the condition
\begin{equation}\label{check2}
\begin{split}
&-\sum_{k=1}^{K-J} \bigg[\sum_{i\in\mathcal{N}_l}(x_i^k-x_i^*)^T(\epsilon_i^k)\bigg]\\
&+\sum_{k=K-J+1}^{K } \bigg[c_u^k(\tilde{x}^*)-c^k_u(x^*)-
\sum_{i\in\mathcal{N}_l}(x_i^k-x_i^*)^T(\epsilon_i^k)\bigg]\leq0,
\end{split}
\end{equation}
holds. Then the average regret of the distribution company   is bounded as
\begin{equation}
\begin{split}
 &\frac{R_u(K)}{K}\leq \frac{1}{K} \bigg[  \frac{1}{\eta_u}P_u+\frac{\eta_u}{2}\sum_{k=1}^{K-J}\Big{\|}\big{(}\nabla c_u^k(x^k)+\epsilon^k\big{)}\Big{\|}_{*}^2\bigg]\\
&+\frac{1}{K} \bigg[ \frac{1}{\eta_u}\tilde{P}_u+\frac{\eta_u}{2}\sum_{k=K-J+1}^{K}\Big{\|}\big{(}\nabla c_u^k(\tilde{x}^k)+\epsilon^k\big{)}\Big{\|}_{*}^2\bigg],
 \end{split}
 \end{equation}
 where
  \begin{equation}
   \tilde{P_u}:= \underset{\tilde{x}\in \tilde{\mathcal{F}}}{\emph{\text{max}}}~L(\tilde{x})-\underset{\tilde{x}\in \tilde{\mathcal{F}}}{\emph{\text{min}}}~L(\tilde{x}).
  \end{equation}
     In  particular, if $\eta=O(1/\sqrt{K})$, then the average regret converges to zero as $K\rightarrow \infty$.
\end{theorem}

The distribution company needs to design the value $J$ and the set $\tilde{\mathcal{F}}_i,~i \in \mathcal{N}_c$ to ensure that the condition (\ref{check2}) is satisfied. The condition (\ref{check2}) imposes requirements on the set of the relaxed feasible charging profiles and the number of days to relax the charging constraints   for the controllable customers.  If the base load remains the same from day to day and the distribution company knows the base load, the value $c_u^k(\tilde{x}^*)$ and $c_u^k(x^*),~k\in \N_{>0}$   can be computed, for example, using the algorithms proposed in \cite{gtl,gwtcl}. The distribution company is able to verify the condition (\ref{check2}) by using the upper bounds of $||x_i^k||$, $||x_i^*||$, and $||\epsilon_i^k||,~i \in \mathcal{N}_l, ~k\in \N_{>0}$ as follows. Given the upper bounds of  $||x_i^k||$, $||x_i^*||$, and $||\epsilon_i^k||,~ i \in \mathcal{N}_l,~ k\in \N_{>0}$, the distribution company requires that the value $J$ and the set $\tilde{\mathcal{F}}_i,~i \in \mathcal{N}_c$ satisfy the inequality 
\begin{equation}
\begin{split}
&\sum_{k=K-J+1}^{K } \bigg(c^k_u(x^*)-c_u^k(\tilde{x}^*)\bigg)\\
& \geq\sum_{k=K-J+1}^{K }
\sum_{i\in\mathcal{N}_l}\bigg(||x_i^k||+||x_i^*||\bigg)||\epsilon_i^k||\\
&+\sum_{k=1 }^{K-J} 
\sum_{i\in\mathcal{N}_l}\bigg(||x_i^k||+||x_i^*||\bigg)||\epsilon_i^k||.
\end{split}
\end{equation}

\section{Numerical Examples}\label{s6}
\begin{figure}[t]
\centering
\includegraphics[width=9cm]{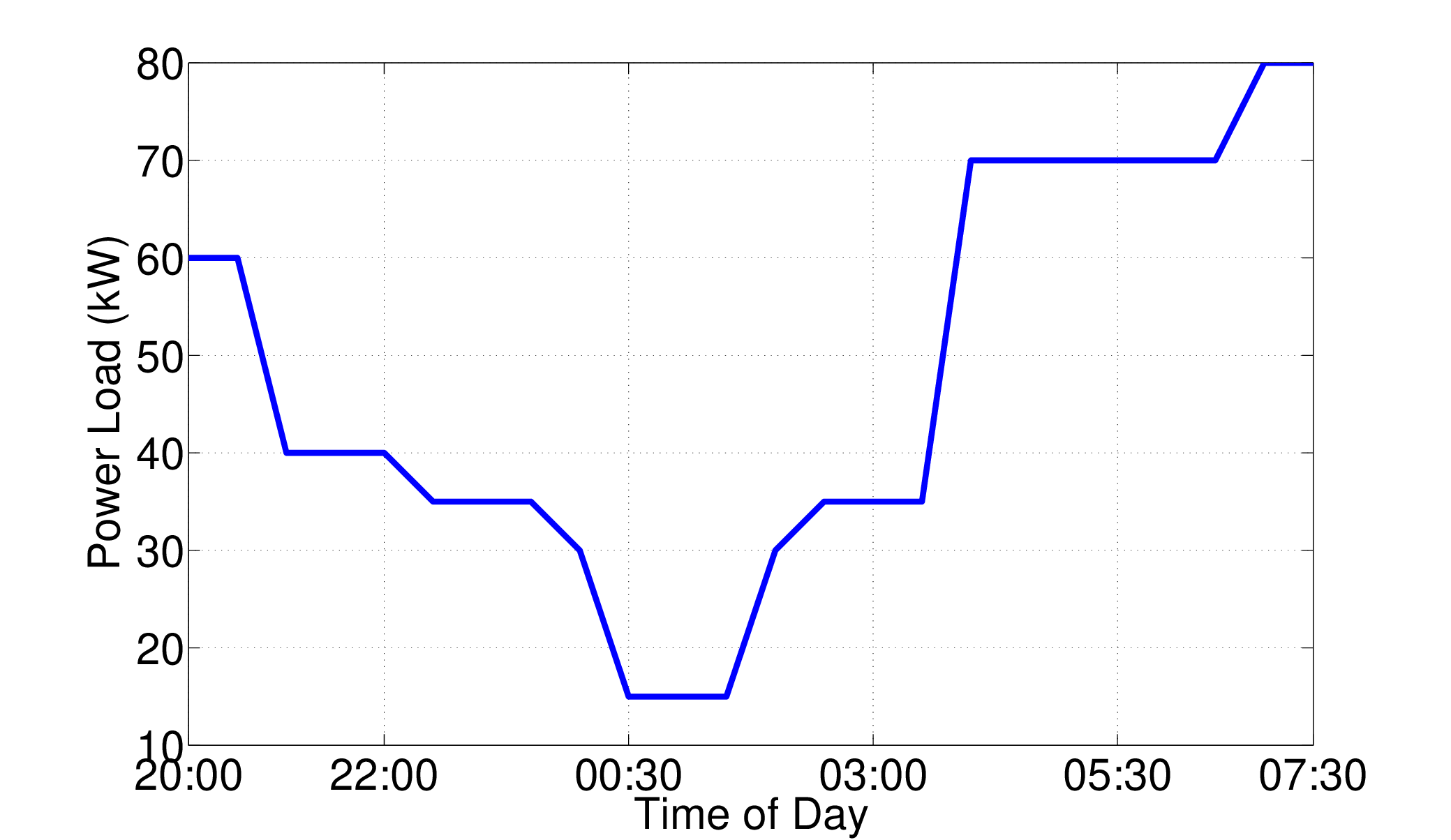}
\caption{Base load profile from 8:00 pm to 7:30 am (the next day).}
\label{0}
\end{figure}

\begin{figure}[t]
\centering
\includegraphics[width=9cm]{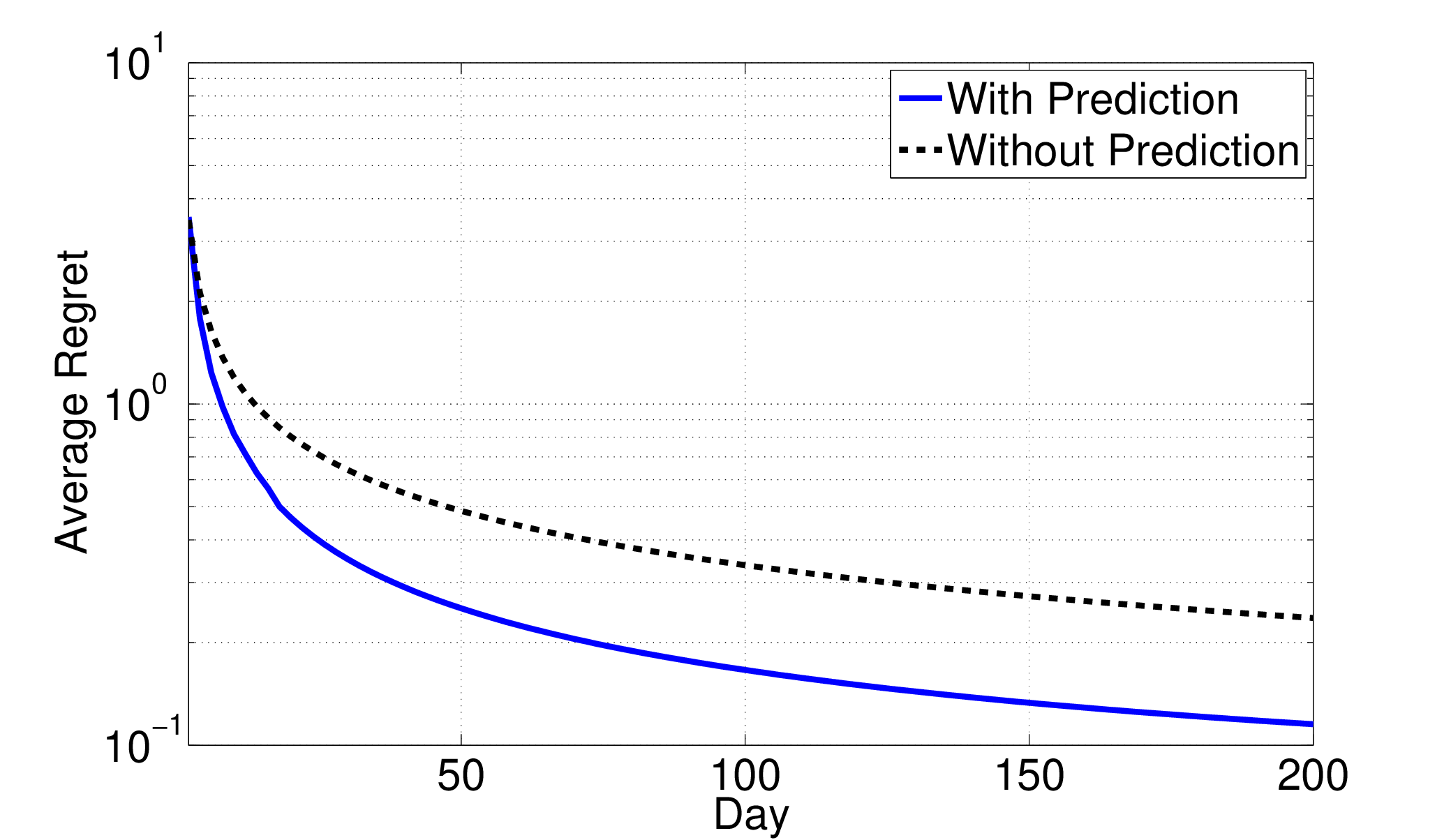}
\caption{Average regrets generated by OMD with and without prediction.}
\label{1}
\vspace{-10mm}
\end{figure}

\begin{figure}[t]
\centering
\includegraphics[width=9cm]{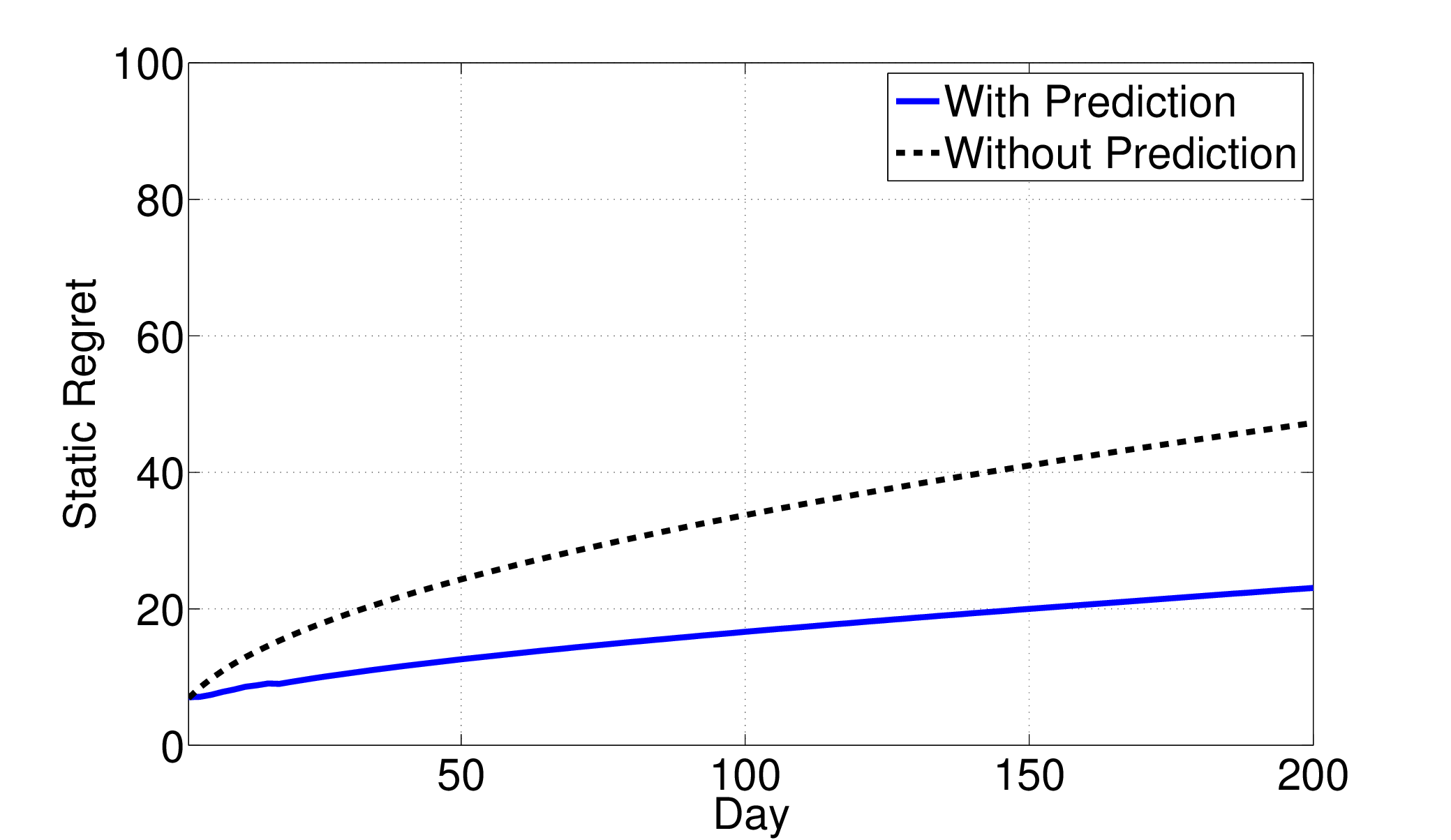}
\caption{Static regrets generated by OMD with and without prediction.}
\label{2}
\end{figure}

Assume that there are 20 customers. A time slot representing an interval of 30 minutes is used. There are $T=24$ time slots. The starting time is set to 8:00 pm. For simplicity, we consider that all EV customers charge their EVs from the $9$th to the $16$th time slots\footnote{The problem formulation (\ref{prob2}) allows EV customers to have different charging constraints. The initial charging time and final charging time mainly depend on the preferences of the EV customers.}. On the first day, the initial charging profiles are assumed to be uniformly distributed over the time slots. The maximum charging rate is set to $x_i^{\mathrm{up}}(t)= 2~ \text{kW},~i\in \mathcal{N}$ and the desired sum $S_i= 10~\text{kW},~i \in \mathcal{N}$. The simulation is carried out for total $K=200$ days. We set the parameters $\eta_i=0.05/\sqrt{K},~i\in \mathcal{N}$. We first examine the convergence of the static regret. The base load profile is given in Figure \ref{0}. Figure \ref{1} shows the trajectories of the regrets with and without the prediction. The prediction  $M_i^k,~i \in \mathcal{N},~k\in \N_{>0}$ is set to $M_i^k(x_i^k)=\frac{1}{k-1}\sum_{\tilde{k}=1,...k-1}\nabla c_i^{\tilde{k}}(x_i^{\tilde{k}}),~i\in \mathcal{N}$. Figure \ref{1} shows that the average regrets converge to zero and the average regret with the prediction converges faster than the one without the prediction. Figure \ref{2} shows the static regrets with and without the prediction. Figure \ref{2} shows that the regrets are sublinear functions of the number of days, which verifies the results in Proposition \ref{prop1}.

\begin{figure}[t]
\centering
\includegraphics[width=9cm]{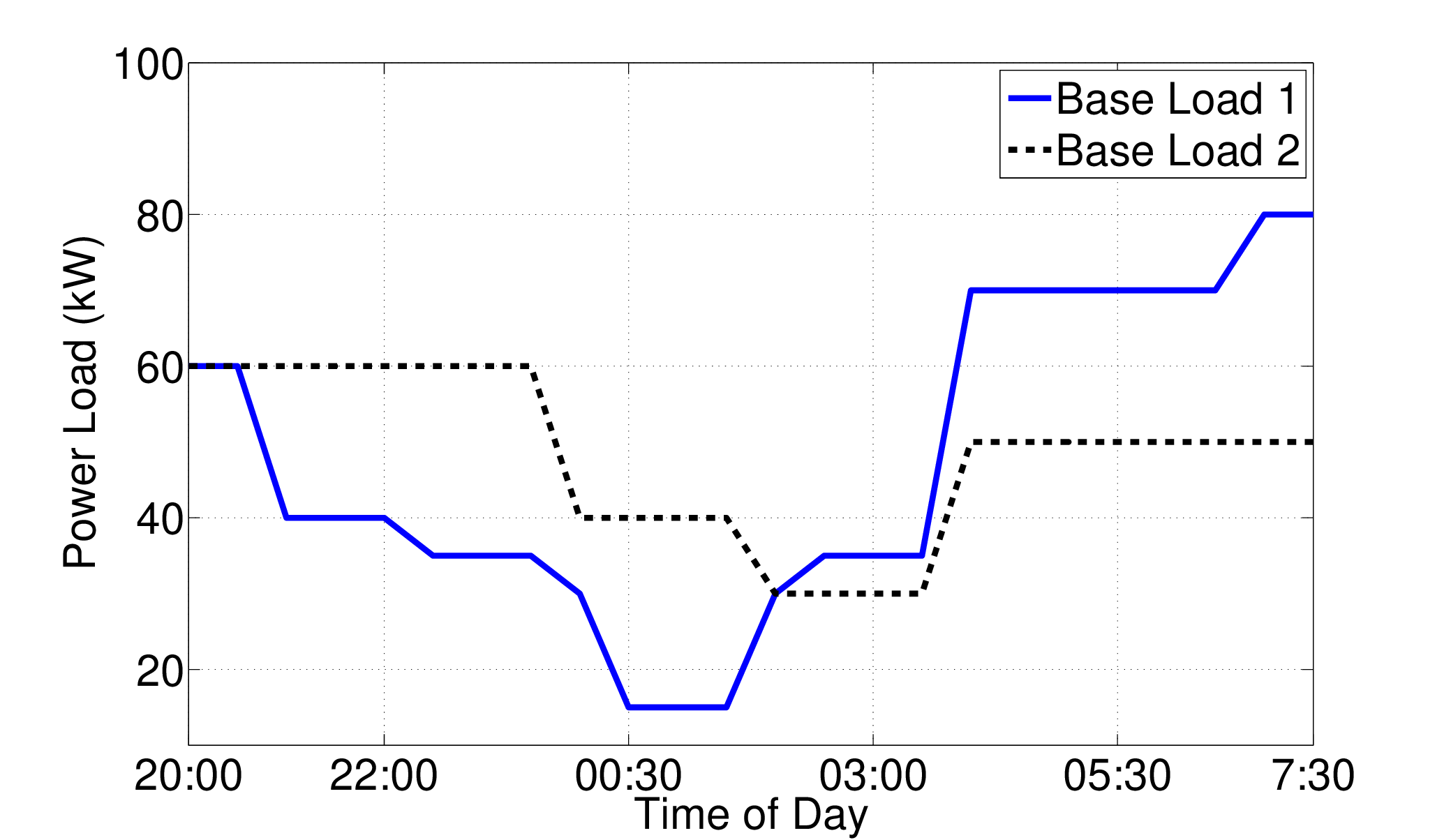}
\caption{Two different base load profiles. The actual base load is realized by switching between the two base load profiles from day to day.}
\label{3}
\vspace{-15mm}
\end{figure}

\begin{figure}[t]
\centering
\includegraphics[width=9cm]{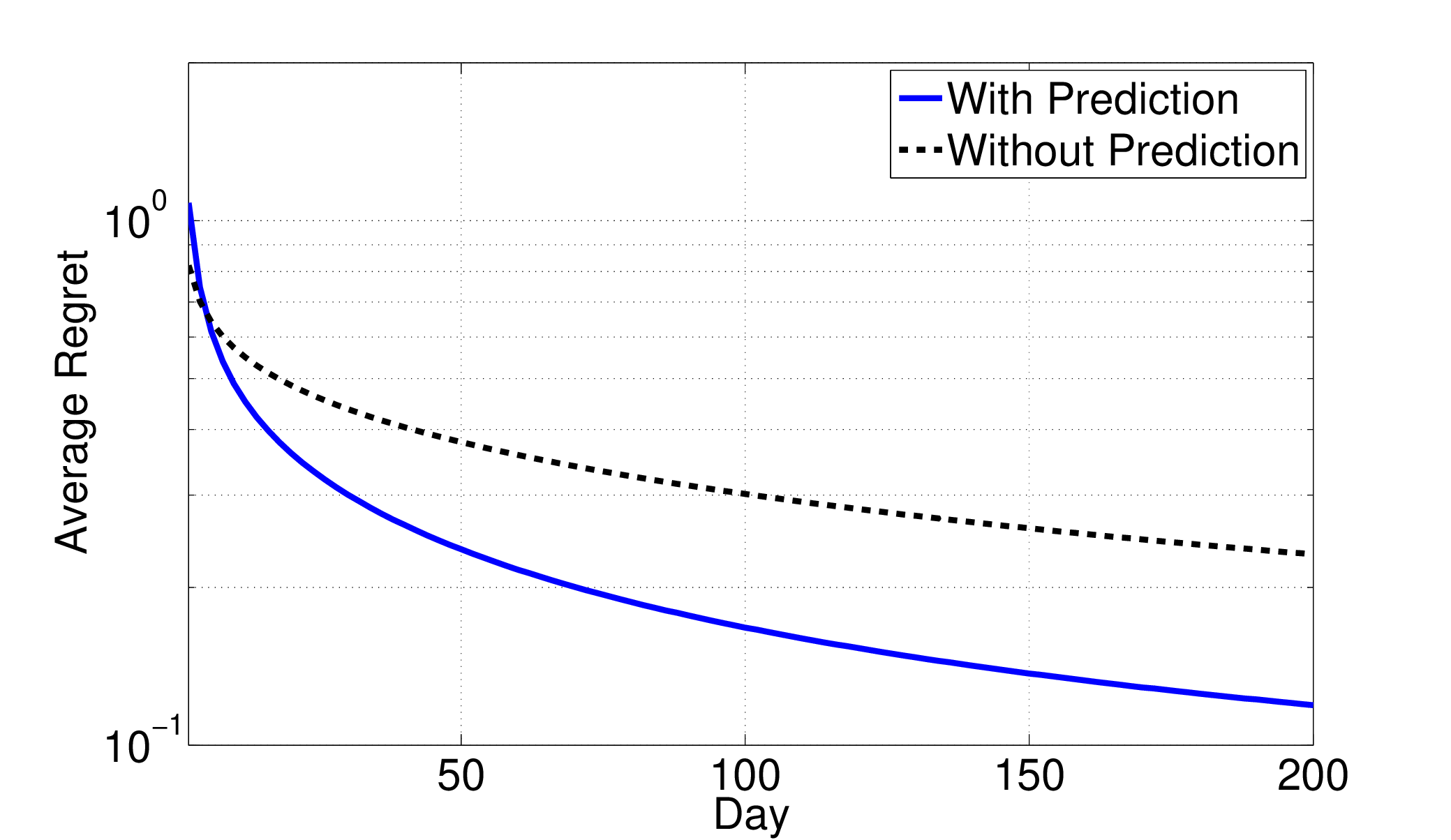}
\caption{Average regrets generated by OMD with and without prediction given the varying base load profile in Figure \ref{3}.}
\label{4}
\end{figure}

\begin{figure}[t]
\centering
\includegraphics[width=9cm]{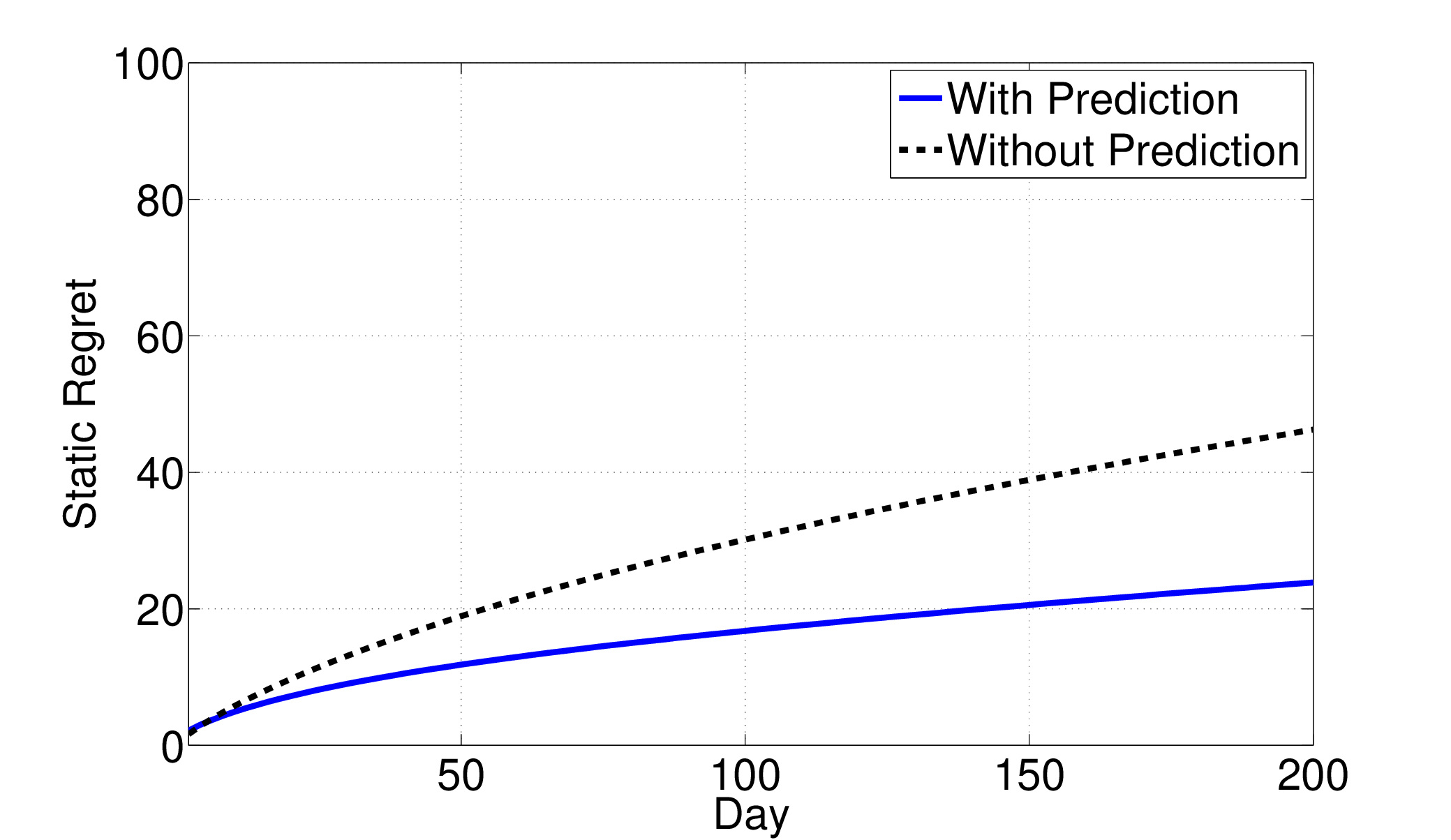}
\caption{Static regrets generated by OMD with and without prediction given the varying base load profile in Figure \ref{3}.}
\label{5}
\end{figure}

We now consider a base load profile which does not remain the same from day to day. The base load is assumed to switch between two base load profiles (see Figure \ref{3}). We set the parameters $\eta_i=0.005/\sqrt{K},~i\in \mathcal{N}$. Figure \ref{4} shows the trajectories of the regrets with and without the prediction given the varying base load in Figure \ref{3}. Figure \ref{4} shows that the average regrets converge to zero and the average regret with the prediction converges faster than the one without the prediction. Figure \ref{5} shows the static regrets with and without the prediction. Figure \ref{5} shows that the regrets are sublinear functions of the number of days. The results in Figure \ref{4} and Figure \ref{5} indicate that despite the fact that the base load is switching and the distribution company is not aware of this varying behavior of the base load, our algorithm still provides updates of the charging profiles having  converging behavior of the regret. Average regret converges to zero means that in the long term, the average performance of the charging schedules $x^k$ generated by the OMD algorithm approaches the performance obtained by $x^*$, where $x^*$ solves (\ref{sub_soln2}).

\begin{figure}[t]
\centering
\includegraphics[width=9cm]{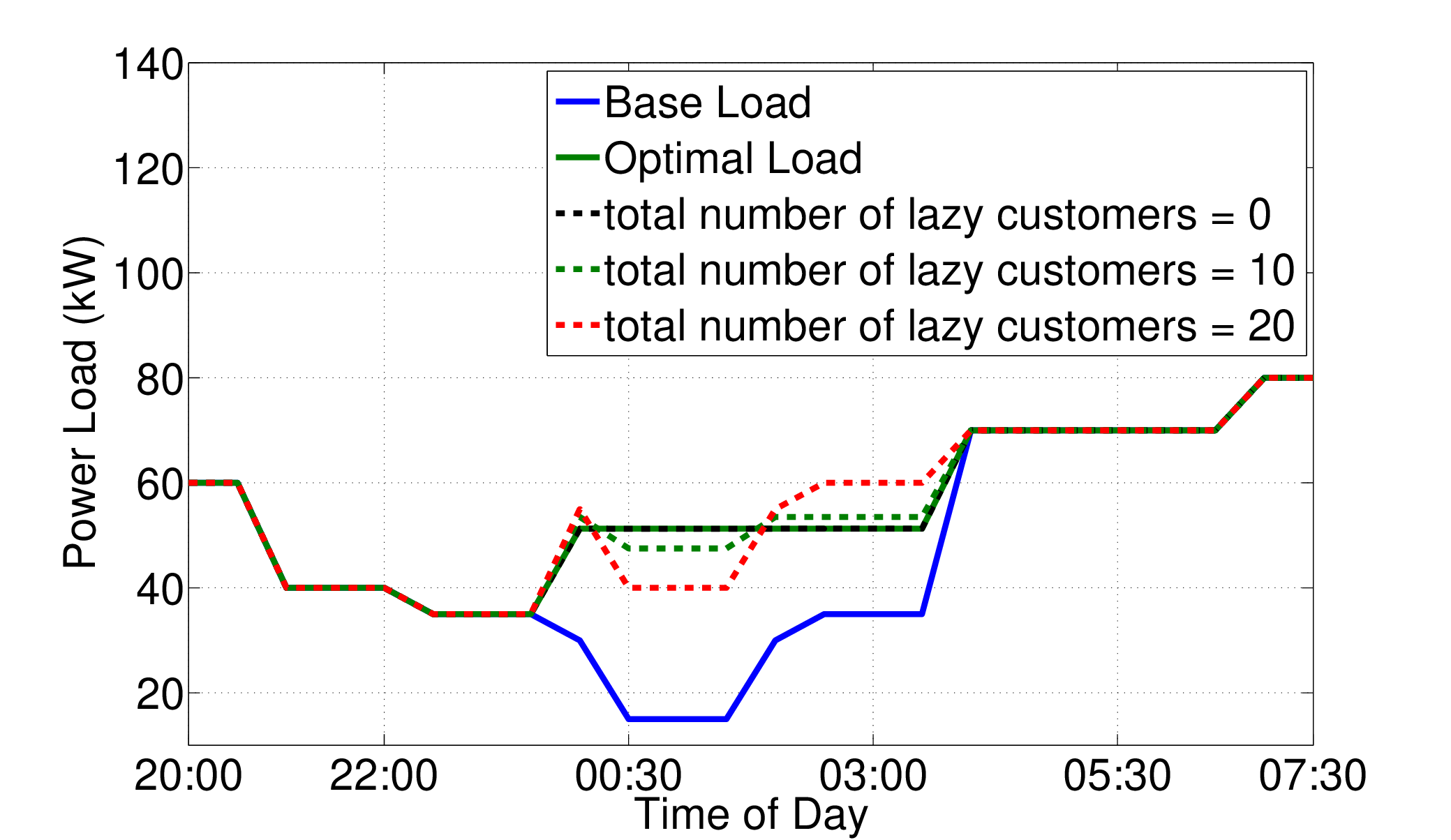}
\caption{The total load profiles with various numbers of inelastic customers, the optimal total load profile, and the base load profile.}
\label{6}
\end{figure}

\begin{figure}[t]
\centering
\includegraphics[width=9cm]{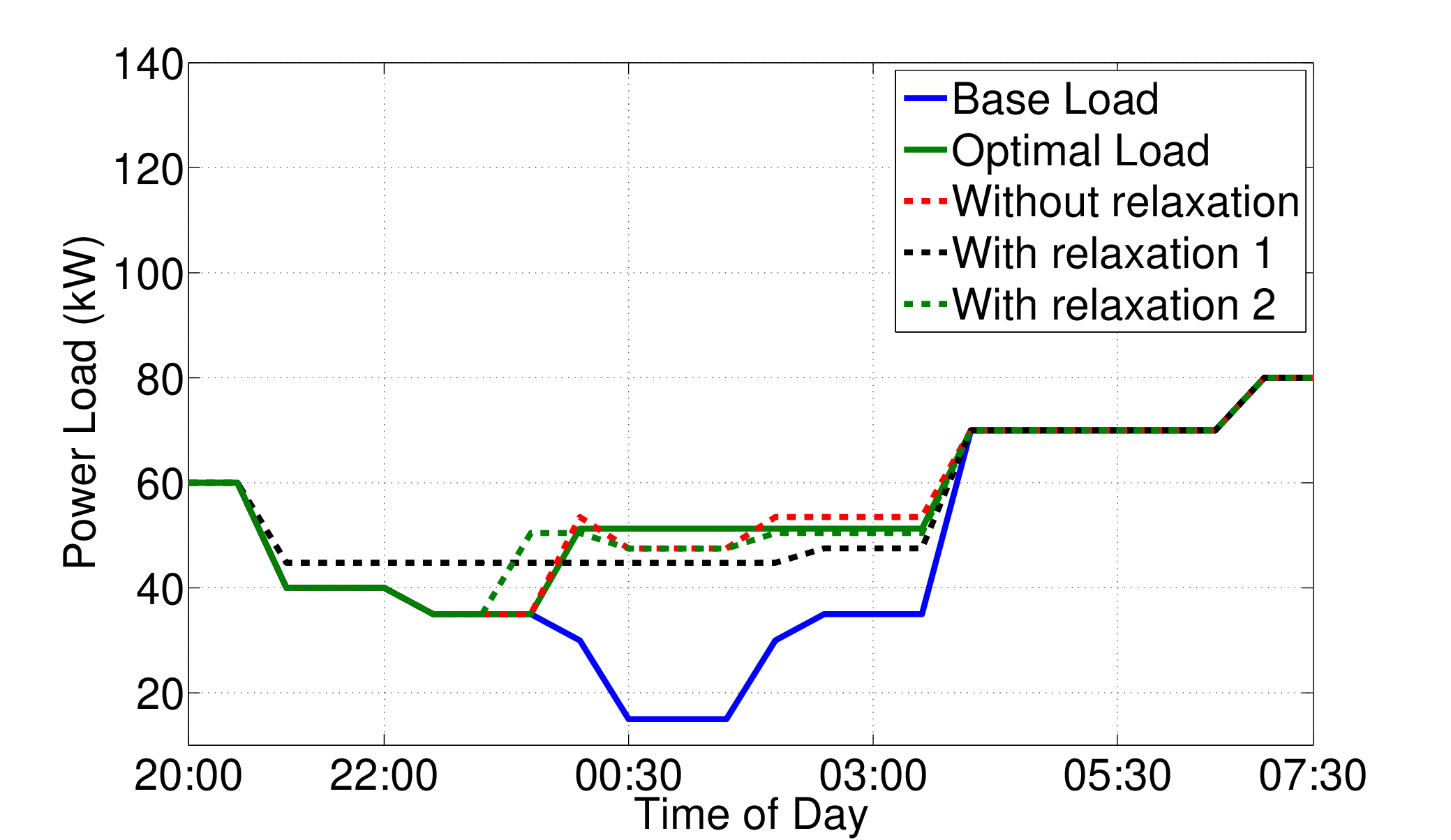}
\caption{The total load profile generated by OMD, the optimal total load profile, base load profile, and the total load profile with the relaxations.}
\vspace{5mm}
\label{7}
\end{figure}

Figure \ref{6} provides the total load profiles at convergence with various number of inelastic customers. Figure \ref{6} shows that as the total number of inelastic customers increases, the variation of the total load profile increases. The result indicates that the inelastic customers perturb the valley-filling load profile. Now suppose that there are 10 inelastic customers and 10 customers with relaxed constraints. We consider two different relaxation strategies. Relaxation 1 represents the strategy that allows the EVs to be charged over the entire time slots (rather than merely between the 9th and the 16th time slots), whereas Relaxation 2 represents the one extending the charging time slots to cover from the 8th to the 17th time slots. Figure \ref{7} shows that as the allowed charging time slots are extended, the total load profile resembles a valley-filing profile, which verifies the effectiveness of the controllable customers to relief stress of the power grid. 

\section{Conclusion}\label{s7}	 
We have designed a framework for distributed charging control of EVs using online learning and online convex optimization. The proposed algorithm can be implemented without low-latency two-way communication between the distribution company and the EV customers, which fits in with the current communication infrastructure and protocols in the smart grid. 
 
\appendices

\section{Proof of Proposition \ref{align}}\label{app_b}

The update (\ref{omd}) yields
\begin{equation}\label{omd31} 
h_i^{k+1}=h_i^k-\eta_i\Bigg{(}\sum_{i}x_i^k+D^k\Bigg{)}, 
\end{equation}
\begin{equation}\label{omd32} 
x_i^{k+1}=\underset{x_i\in \mathcal{F}_i}{\text{argmin}}~ \eta_i x_i^TM_i^{k+1}+\|x_i-h_i^{k+1}\|^2, 
\end{equation}
where $M_i^k$ is the prediction of the value $\sum_{i}x_i^k+D^k$.

 The update (\ref{omdu}) yields
\begin{equation}\label{propcom1}
  \begin{bmatrix}
       h_1^{k+1}\\
       h_2^{k+1}\\ 
      \vdots\\
       h_N^{k+1}
    \end{bmatrix}_{NT}=
   \begin{bmatrix}
       h_1^{k}\\
       h_2^{k}\\ 
      \vdots\\
       h_N^{k}
    \end{bmatrix}_{NT}
    -\eta_u
    \begin{bmatrix}
      2(D^k+\sum_{i\in \mathcal{N}}x_i^k)\\
      2(D^k+\sum_{i\in \mathcal{N}}x_i^k)\\ 
      \vdots\\
      2(D^k+\sum_{i\in \mathcal{N}}x_i^k)
    \end{bmatrix}_{NT},\\\\
\end{equation} 
\begin{equation}\label{propcom2}
      x^{k+1}=\underset{x\in \mathcal{F}}{\text{argmin}}~ \eta_u x^TM_u^{k+1}+\|x-h_u^{k+1}\|^2,
    \end{equation} 
where $h_u^{k+1}=[ h_1^{k+1},..., h_N^{k+1}]^T$ and $M_u^k$ is the prediction of

\begin{equation}
\begin{bmatrix}
      2(D^k+\sum_{i\in \mathcal{N}}x_i^k)\\
      2(D^k+\sum_{i\in \mathcal{N}}x_i^k)\\ 
      \vdots\\
      2(D^k+\sum_{i\in \mathcal{N}}x_i^k)
       \end{bmatrix}_{NT.}
 \end{equation}
 
By comparing (\ref{propcom1}), (\ref{propcom2}) with (\ref{omd31}),  (\ref{omd32}) and substituting $\eta_u=\frac{1}{2}\eta_i,~i\in \mathcal{N}$, we have that the updates (\ref{propcom2})  and   (\ref{omd32}) are identical. Since the updates coincide, the average regret of the distribution company converges to zero as $k\rightarrow \infty$.

\section{Proof of Theorem \ref{thm:tracking}}\label{app_d}
The proof technique that we use to derive the tracking regret bound in (\ref{tracking_upbd}) is similar to the one used in \cite[Theorem 4]{hw}. The main step is to bound the difference $c_u^k(x^k)-c_u^k(x^{k*})$ instead of the difference $c_u^k(x^k)-c_u^k(x^{*})$ that is considered in the static regret (\ref{upbd2}). However, in \cite{hw}, the authors derive the tracking regret bounds for a different regret minimization algorithm rather than OMD. 

Following the proof of \cite[Lemma 2]{rs}, we have
\begin{equation}
c_u^k(x^k)-c_u^k(x^{k*})\leq(x^k-x^{k*})^T\nabla c_u^k(x^k),
\end{equation}
and 
\begin{equation}\label{thm42_update}
\begin{split}
&(x^k-x^{k*})^T\nabla c_u^k(x^k)\leq \frac{\eta_u}{2} \|\nabla c_u^k(x^k)-M_u^k\|_{*}^2\\
&+\frac{1}{\eta_u}\big{(}D_{L_u}(x^{k*},h_u^{k})-D_{L_u}(x^{k*},h_u^{k+1})\big{)}.
\end{split}
\end{equation}
Furthermore, 
\begin{equation}
\begin{split}
&D_{L_u}(x^{k*},h_u^{k})-D_{L_u}(x^{k*},h_u^{k+1})\\
&=L_u(x^{k*})-L_u(h_u^{k})-\nabla L_u(h_u^{k})^T(x^{k*}-h_u^k)\\
&-L_u(x^{k*})+L_u(h_u^{k+1})+\nabla L_u(h_u^{k+1})^T(x^{k*}-h_u^{k+1})\\
&=L_u(h_u^{k+1})-L_u(h_u^{k})+\nabla L_u(h_u^{k+1})^T(x^{k+1*}-h_u^{k+1})\\
&-\nabla L_u(h_u^{k})^T(x^{k*}-h_u^{k})-\nabla L_u(h_u^{k+1})^T(x^{k+1*}-x^{k*}).
\end{split}
\end{equation}
The remainder of the proof is followed by summing over $k=1,...,K$ and collecting terms.

\section{Proof of Theorem \ref{lazy}}\label{app_e}
Following the proof of \cite[Lemma 2]{rs}, we have
\begin{equation}
c_u^k(x^k)-c_u^k(x^*)\leq(x^k-x^{*})^T\nabla c_u^k(x^k),
\end{equation}
and 
\begin{equation}\label{thm42_update}
\begin{split}
&(x^k-x^{*})^T\nabla c_u^k(x^k)\leq \frac{\eta_u}{2} \|\nabla c_u^k(x^k)\|_{*}^2\\
&+\frac{1}{\eta_u}\big{(}D_{L_u}(x^{*},h_u^{k})-D_{L_u}(x^{*},h_u^{k+1})\big{)},
\end{split}
\end{equation}
where 
\begin{equation} 
\nabla c_u^k(x^k)=
\begin{bmatrix}
      2(D^k+\sum_{i\in \mathcal{N}}x_i^k)\\ 
      \vdots\\
      2(D^k+\sum_{i\in \mathcal{N}}x_i^k)
       \end{bmatrix}_{NT}.
        \end{equation} 
For each $i\in \mathcal{N}_l$, the cost function $c_i^k$ is selected as a constant function noted in Section \ref{4a}. The corresponding gradient of the lazy customer's cost function is zero, namely, for lazy customer $i\in \mathcal{N}_l$,  
\begin{equation}
\nabla c_i^k(x^k)=
      \Bigg{(}D^k+\sum_{i\in \mathcal{N}}x_i^k\Bigg{)}+\epsilon_i^k=0,
\end{equation}
where $\epsilon_i^k=-( D^k+\sum_{i\in \mathcal{N}}x_i^k)$.
Move $\epsilon_i^k,~i\in \mathcal{N}_l $ to the right hand side of the inequality in (\ref{thm42_update}). The remainder of the proof is followed by summing over $k=1,...,K$ and collecting terms.

\section{Proof of Theorem \ref{lazy2}}\label{app_f}
First observe that, for $i\in \mathcal{N}_l$,
\begin{equation}
c_u^k(\tilde{x}^k)-c_u^k(x^*)=[c_u^k(\tilde{x}^k)-c_u^k(\tilde{x}^*)]+[c_u^k(\tilde{x}^*)-c_u^k(x^*)].
\end{equation}

Following the proof of \cite[Lemma 2]{rs},  
\begin{equation} 
\begin{split}
&c_u^k(\tilde{x}^k)-c_u^k(\tilde{x}^*)\leq (\tilde{x}^k-\tilde{x}^{*})^T\nabla c_u^k(\tilde{x}^k)\\
&\leq \frac{\eta_u}{2} \|\nabla c_u^k(\tilde{x}^k)\|_{*}^2\\
&+\frac{1}{\eta_u}\big{(}D_{L_u}(\tilde{x}^{*},h_u^{k})-D_{L_u}(\tilde{x}^{*},h_u^{k+1})\big{)}.
\end{split}
\end{equation}

The remainder of the proof is followed by summing over $k=1,...,K$ and collecting terms.

\end{document}